# A note on the transformation of the linear differential equation into a system of the first order equations

M.I. Ayzatsky[1]

National Science Center Kharkov Institute of Physics and Technology (NSC KIPT), 610108, Kharkov, Ukraine

A generalization of the already studied transformations of the linear differential equation into a system of the first order equations is given. The proposed transformation gives possibility to get new forms of the $N$-dimensional system of first order equations that can be useful for analysis of the solutions of the $N$-th order differential equations. In particular, for the third-order linear equation the nonlinear second-order equation that plays the same role as the Riccati equation for second-order linear equation is obtained.

## 1 Introduction

"It is sometimes possible to transform a given differential equation into an equation which can be integrated directly. In many cases the new equation is derived from the original one by a change of variables, sometimes it is derived from it by differentiation. It is generally possible to use a solution of the new equation to obtain a solution of the original equation. If, however, the new equation cannot be solved in a simple manner by means of known functions, we must not conclude that the transformation is of no value, because a differential equation is frequently used to define a new functions, and then it is important to know which types of differential equations can be solved with the aid of the new functions. Transformations which change a given differential equation into another differential equation of the same general type are of some interest, as they frequently indicate properties of the functions defined by the given differential equations." ([1],p.60)

It is well known that the $N$-th order differential equation[2]

$$y^{(N)} + f_{N-1}(t) y^{(N-1)} + \ldots \ldots f_1(t) y' + f_0(t) y + f(t) = 0 \tag{1}$$

can be converted into an $N$-dimensional system of first order equations. There are various reasons for doing this, one being that a first order system is much easier to solve numerically as the most differential equations we encounter in physics, economics and engineering don't have exact solutions.

The most common method is introducing a number of new variables

$$y^{(m)} = y_m, \quad m = 1, \ldots N-1 \tag{2}$$

and representing (1) as

$$Y' = MY + H, \tag{3}$$

where

$$H = \left(0, 0 \ldots, -f(t)\right)^T,$$
$$Y = \left(y, y_1, \ldots, y_{N-1}\right)^T, \tag{4}$$

---

[1] M.I. Aizatskyi, N.I.Aizatsky; aizatsky@kipt.kharkov.ua
[2] A prime denotes differentiation with respect to t

$$M = \begin{pmatrix} 0 & 1 & 0 & \ldots & 0 \\ 0 & 0 & 1 & \ldots & 0 \\ \ldots & \ldots & \ldots & \ldots & \ldots \\ 0 & 0 & 0 & \ldots & 1 \\ f_0 & f_1 & f_2 & \ldots & f_{N-1} \end{pmatrix}. \quad (5)$$

This system of equations can be subjected to further transformations (see, for example, [2,3,4,5,6]).

For the second order differential equation
$$y'' + f_1(t)y' + f_0(t)y = 0. \quad (6)$$
there is a set of transformation that reduce this equation to the form ([7])
$$y' = P(t)y + Q(t)z, \quad z' = R(t)y + S(t)z, \quad (7)$$
where $z(t)$ is an additional unknown function. In the simplest case $P(t) \equiv 0$, $Q(t) \equiv 1$ we obtain the described above transformation (3).

All these transformations assume that along with the unknown function $y(t)$, $(N-1)$ additional functions are introduced.

There is another kind of transformation that consists in presentation of the solution $y(t)$ of the equation (1) or its derivatives as the sum of the $N$ new unknown functions multiplied by the known functions[3]. By introducing $N$ new unknowns $y_n(t)$ instead of the one $y(t)$, we can impose $(N-1)$ additional conditions. Such approach can give new form of the $N$-dimensional system of first order equations, equivalent to the equation (1). For the case $N=2$ it is widely used in a theory of electromagnetic waves in stratified media (see, for example, [2,8,9,10]). For the case $N=2$ there are some other transformations of the equation (6) (see, for example, [5]).

In this note we propose some generalization of the already studied transformations. Analysis of literature [1,5,7,11,12,13,14,15,16,17,18,19,20] shows that it apparently has not been described earlier. Proposed approach is similar to the transformation a second-order linear difference equation by using two arbitrary sequences of numbers [21].

Bellow, for simplicity and clarity, we consider the case $N=2$ and $N=3$. More general approach for arbitrary $N$ is given in Appendix.

## 2 Linear differential equation of the second order

Suppose the solution of the equation
$$y'' + f_1(t)y' + f_0(t)y + f(t) = 0 \quad (8)$$
can be represented as the sum of two new functions
$$y(t) = y_1(t) + y_2(t). \quad (9)$$

Additional condition we write in the form
$$y' = g_1(t)y_1(t) + g_2(t)y_2(t), \quad (10)$$
where $g_1(t)$ and $g_2(t)$ are arbitrary continuous functions having continuous derivatives.

For given $g_1, g_2$ ($g_1 \neq g_2$) this representation is unique. Indeed, from (9) and (10) we obtain
$$y_1 = \frac{1}{g_2 - g_1}\{yg_2 - y'\},$$
$$y_2 = -\frac{1}{g_2 - g_1}\{yg_1 - y'\}. \quad (11)$$

For the Cauchy problem these formulae gives initial values for $y_1(t_1)$ and $y_2(t_1)$.

---

[3] Usually, they are the WKB functions

Substitutions of (9) into (10) and (9),(10) into (8) give

$$g_1 y'_1 + g_2 y'_2 = -g'_1 y_1 - g'_2 y_2 - f_1(g_1 y_1 + g_2 y_2) - f_0(y_1 + y_2) - f,$$
$$y'_1 + y'_2 = (g_1 y_1 + g_2 y_2). \qquad (12)$$

As $g_1(t) \neq g_2(t)$, from this system we can find the derivatives of $y_1$ and $y_2$

$$(g_1 - g_2) y'_1 = -y_1(g'_1 + f_1 g_1 + f_0 + g_2 g_1) - y_2(g'_2 + f_1 g_2 + f_0 + g_2^2) - f,$$
$$(g_1 - g_2) y'_2 = y_1(g'_1 + f_1 g_1 + f_0 + g_1^2) + y_2(g'_2 + f_1 g_2 + f_0 + g_2 g_1) + f. \qquad (13)$$

This system is the basis of the proposed transformation in the case of the second-order linear differential equation (8).

We would like to emphasize that the functions $g_1(t)$ and $g_2(t)$ are the arbitrary ones, and we do not impose a condition that the functions $y_1(t)$ and $y_2(t)$ separately are the solutions of equation(8), but their sum must be the solution of this equation.

Introducing new functions

$$x_1 = g'_1 + f_1 g_1 + f_0 + g_1^2,$$
$$x_2 = g'_2 + f_1 g_2 + f_0 + g_2^2, \qquad (14)$$

equations (13) can be rewritten in the form

$$y'_1 = y_1 g_1 - y_1 \frac{x_1}{(g_1 - g_2)} - y_2 \frac{x_2}{(g_1 - g_2)} - \frac{f}{(g_1 - g_2)},$$
$$y'_2 = y_2 g_2 + y_1 \frac{x_1}{(g_1 - g_2)} + y_2 \frac{x_2}{(g_1 - g_2)} + \frac{f}{(g_1 - g_2)}. \qquad (15)$$

Let's consider several cases when special choosing of $g_1(t)$ and $g_2(t)$ gives useful results.

If we choose

$$g_{1,2} = \rho_{1,2}(t). \qquad (16)$$

where $\rho_{1,2}$ are the solutions of the characteristic equation

$$\rho^2 + f_1(t)\rho + f_0(t) = 0, \qquad (17)$$

then the system (15) takes the form

$$y'_1 = y_1 \rho_1 - y_1 \frac{\rho'_1}{(\rho_1 - \rho_2)} - y_2 \frac{\rho'_2}{(\rho_1 - \rho_2)} - \frac{f}{(\rho_1 - \rho_2)},$$
$$y'_2 = y_2 \rho_2 + y_1 \frac{\rho'_1}{(\rho_1 - \rho_2)} + y_2 \frac{\rho'_2}{(\rho_1 - \rho_2)} + \frac{f}{(\rho_1 - \rho_2)}. \qquad (18)$$

For the case $f_1(t) \equiv 0$, $f(t) \equiv 0$

$$y'_1 = \left( i\sqrt{f_0} - \frac{f'_0}{4 f_0} \right) y_1 + y_2 \frac{f'_0}{4 f_0},$$
$$y'_2 = -\left( i\sqrt{f_0} + \frac{f'_0}{4 f_0} \right) y_2 + y_1 \frac{f'_0}{4 f_0}. \qquad (19)$$

This is the well known system that is a basis for coupled wave model and WKB approximation [3,4,22,23,24,25,26]. We would like to note again that $y_1(t)$ and $y_2(t)$ separately are not the solutions of equation (8).

Neglecting the second terms on the right-hand side of the equations(19), we obtain the well known WKB result

$$y_{1,2} = y_{1,2}(t_1)\left(\frac{f_0(t_1)}{f_0(t)}\right)^{1/4}\exp\left(\pm i\int_{t_1}^{t}\sqrt{f_0}dt'\right). \quad (20)$$

If $g_1(t)$ and $g_2(t)$ are the two different solutions of the Riccati equation [27,28]

$$g'_{1,2} + f_1 g_{1,2} + f_0 + g_{1,2}^2 = 0, \quad (21)$$

the system (13) transforms into

$$y'_1 = g_1 y_1 - \frac{f}{g_1 - g_2},$$
$$y'_2 = g_2 y_2 + \frac{f}{g_1 - g_2}. \quad (22)$$

For the homogeneous equation ($f(t) \equiv 0$)

$$y_{1,2} = y_{1,2}(t_1)\exp\left(\int_{t_1}^{t} g_{1,2}dt'\right). \quad (23)$$

The derivative of the Wronskian of these functions equals

$$W' = y_1(t_1)y_2(t_1)\left[(g_2 - g_1)' + (g_2^2 - g_1^2)\right]\exp\left(\int_{t_1}^{t}(g_1 + g_2)dt'\right). \quad (24)$$

Inserting the derivatives of $g_{1,2}$ from (21) into this expression, we obtain

$$W' = -f_1(t)W, \quad (25)$$

that proves the linearly independence of the functions $y_{1,2}$.

From the Riccati equations (21) we can find that ($g_1(t) \neq g_2(t)$)

$$\frac{1}{g_1 - g_2}(g_1 - g_2)' = -(g_1 + g_2) - f_1. \quad (26)$$

Denoting the difference of functions $g_1$ and $g_2$ as

$$(g_1 - g_2) = 2iq, \quad (27)$$

then from (26) and (27) we obtain.

$$g_{1,2} = \pm iq - \frac{1}{2q}q' - \frac{f_1}{2}. \quad (28)$$

The solutions (23) can be rewritten in the form that is commonly used in the phase-integral method [29,30]

$$y_{1,2} = y_{1,2}(t_1)\left(\frac{q(t_1)}{q(t)}\right)^{1/2}\exp\left(\pm i\int_{t_1}^{t}qdt' - \frac{1}{2}\int_{t_1}^{t}f_1 dt'\right). \quad (29)$$

Inserting (28) into the Riccati equations (21), we find the equation for the $q$ function

$$g'_{1,2} + f_1 g_{1,2} + f_0 + g_{1,2}^2 = -\frac{1}{2q}q'' + \frac{3}{4q^2}(q')^2 + f_0 - q^2 - \frac{f'_1}{2} - \frac{f_1^2}{4} = 0, \quad (30)$$

which for $f_1(t) \equiv 0$ takes the well known form [30]

$$q^{-3/2}(q^{-1/2})'' + \frac{f_0}{q^2} - 1 = 0. \quad (31)$$

If $g_1(t)$ and $g_2(t)$ are the solutions of the system of nonlinear differential equations

$$g'_1 + f_1 g_1 + f_0 + g_2 g_1 = 0,$$
$$g'_2 + f_1 g_2 + f_0 + g_2 g_1 = 0, \quad (32)$$

then the system (13) transforms into the "strong" coupling system

$$y_1' = y_2 g_2 - \frac{f}{(g_1 - g_2)},$$
$$y_2' = y_1 g_1 + \frac{f}{(g_1 - g_2)}. \tag{33}$$

The system (32) can be transformed into the Riccati equation

$$g_1' + \left[ f_1 - C \exp\left( -\int_{t_t}^{t} f_1 dt' \right) \right] g_1 + f_0 + g_1^2 = 0 \tag{34}$$

and the relationship between $g_1$ and $g_2$

$$g_1 - g_2 = C \exp\left( -\int_{t_t}^{t} f_1 dt' \right). \tag{35}$$

$C \neq 0$ is the arbitrary constant.

## 2 Linear differential equation of the third-order

Let's represent the solution of the linear third-order equation (results of detail study of its solutions see in [31])

$$y''' + f_2(t) y'' + f_1(t) y' + f_0(t) y + f(t) = 0. \tag{36}$$

as the sum of three new functions

$$y(t) = y_1(t) + y_2(t) + y_3(t). \tag{37}$$

Additional conditions we write in the form

$$y' = g_{1,1}(t) y_1(t) + g_{1,2}(t) y_2(t) + g_{1,3}(t) y_3(t),$$
$$y'' = g_{2,1}(t) y_1(t) + g_{2,2}(t) y_2(t) + g_{2,3}(t) y_3(t), \tag{38}$$

where $g_{m,n}(t)$ ($m=1,2;\ n=1,2,3$) are arbitrary continuous functions with continuous derivatives and

$$D = \begin{vmatrix} 1 & 1 & 1 \\ g_{1,1} & g_{1,2} & g_{1,3} \\ g_{2,1} & g_{2,2} & g_{2,3} \end{vmatrix} \neq 0. \tag{39}$$

Making the transformations that are given in Appendix, we obtain a system of the first-order linear differential equations

$$Dy_1' = y_1 \left[ g_{1,1} D + (g_{2,2} - g_{2,3})(x_1 - (g_{1,1})^2) - (g_{1,3} - g_{1,2})(x_4 + g_{1,1} g_{2,1}) \right] +$$
$$+ y_2 \left[ (g_{2,2} - g_{2,3})(x_2 - (g_{1,2})^2) - (g_{1,3} - g_{1,2})(x_5 + g_{1,2} g_{2,2}) \right] + \tag{40}$$
$$+ y_3 \left[ (g_{2,2} - g_{2,3})(x_3 - (g_{1,3})^2) - (g_{1,3} - g_{1,2})(x_6 + g_{1,3} g_{2,3}) \right] - c(g_{1,3} - g_{1,2}),$$

$$Dy_2' = y_1 \left[ (g_{2,3} - g_{2,1})(x_1 - (g_{1,1})^2) - (g_{1,1} - g_{1,3})(x_4 + g_{1,1} g_{2,1}) \right] +$$
$$+ y_2 \left[ g_{1,2} D + (g_{2,3} - g_{2,1})(x_2 - (g_{1,2})^2) - (g_{1,1} - g_{1,3})(x_5 + g_{1,2} g_{2,2}) \right] + \tag{41}$$
$$+ y_3 \left[ (g_{2,3} - g_{2,1})(x_3 - (g_{1,3})^2) - (g_{1,1} - g_{1,3})(x_6 + g_{1,3} g_{2,3}) \right] - c(g_{1,1} - g_{1,3}),$$

$$Dy'_3 = y_1\left[(g_{2,1} - g_{2,2})(x_1 - (g_{1,1})^2) - (g_{1,2} - g_{1,1})(x_4 + g_{1,1}g_{2,1})\right] +$$
$$+ y_2\left[(g_{2,1} - g_{2,2})(x_2 - (g_{1,2})^2) - (g_{1,2} - g_{1,1})(x_5 + g_{1,2}g_{2,2})\right] + \quad (42)$$
$$+ y_3\left[g_{1,3}D + (g_{2,1} - g_{2,2})(x_3 - (g_{1,3})^2) - (g_{1,2} - g_{1,1})(x_6 + g_{1,3}g_{2,3})\right] - c(g_{1,2} - g_{1,1}),$$

where the following notations were introduced

$$\begin{aligned} x_1 &= (g_{2,1} - g'_{1,1}), \\ x_2 &= (g_{2,2} - g'_{1,2}), \\ x_3 &= (g_{2,3} - g'_{1,3}), \\ x_4 &= (g'_{2,1} + f_2 g_{2,1} + f_1 g_{1,1} + f_0), \\ x_5 &= (g'_{2,2} + f_2 g_{2,2} + f_1 g_{1,2} + f_0), \\ x_6 &= (g'_{2,3} + f_2 g_{2,3} + f_1 g_{1,3} + f_0). \end{aligned} \quad (43)$$

If we choose

$$\begin{aligned} g_{1,n} &= \rho_n, \\ g_{2,n} &= \rho_n^2, \ n = 1,2,3 \end{aligned} \quad (44)$$

where $\rho_n$ are the solutions of the equation

$$\rho^3 + f_2(t)\rho^2 + f_1(t)\rho + f_0(t) = 0, \quad (45)$$

then the system (42) takes the form

$$\begin{aligned} y'_1 &= y_1\rho_1 - y_1\rho'_1\left(\frac{1}{\rho_1 - \rho_3} + \frac{1}{\rho_1 - \rho_2}\right) - y_2\rho'_2\left(\frac{1}{\rho_1 - \rho_2} + \frac{1}{\rho_3 - \rho_1}\right) - \\ &\quad - y_3\rho'_3\left(\frac{1}{\rho_1 - \rho_3} + \frac{1}{\rho_2 - \rho_1}\right) - f\frac{\rho_3 - \rho_2}{D}, \\ y'_2 &= y_2\rho_2 - y_1\rho'_1\left(\frac{1}{\rho_2 - \rho_2} + \frac{1}{\rho_3 - \rho_2}\right) - y_2\rho'_2\left(\frac{1}{\rho_2 - \rho_1} + \frac{1}{\rho_2 - \rho_3}\right) - \\ &\quad - y_3\rho'_3\left(\frac{1}{\rho_1 - \rho_3} + \frac{1}{\rho_2 - \rho_1}\right) - f\frac{\rho_1 - \rho_3}{D}, \\ y'_3 &= y_3\rho_3 - y_1\rho'_1\left(\frac{1}{\rho_3 - \rho_3} + \frac{1}{\rho_2 - \rho_3}\right) - y_2\rho'_2\left(\frac{1}{\rho_3 - \rho_2} + \frac{1}{\rho_1 - \rho_3}\right) - \\ &\quad - y_3\rho'_3\left(\frac{1}{\rho_3 - \rho_1} + \frac{1}{\rho_3 - \rho_3}\right) - f\frac{\rho_2 - \rho_1}{D}, \end{aligned} \quad (46)$$

where $D = (\rho_3 - \rho_1)(\rho_3 - \rho_2)(\rho_2 - \rho_1)$. This system coincides with the one that was obtained in [5] by another method. In the WKB approximation

$$\begin{aligned} y'_1 &= y_1\rho_1 - y_1\rho'_1\left(\frac{1}{\rho_1 - \rho_3} + \frac{1}{\rho_1 - \rho_2}\right) - f\frac{\rho_3 - \rho_2}{D}, \\ y'_2 &= y_2\rho_2 - y_2\rho'_2\left(\frac{1}{\rho_2 - \rho_1} + \frac{1}{\rho_2 - \rho_3}\right) - f\frac{\rho_1 - \rho_3}{D}, \\ y'_3 &= y_3\rho_3 - y_3\rho'_3\left(\frac{1}{\rho_3 - \rho_1} + \frac{1}{\rho_3 - \rho_3}\right) - f\frac{\rho_2 - \rho_1}{D}, \end{aligned} \quad (47)$$

it also coincides with results of asymptotic analysis (see, for example, [5])

If we choose the functions $g_{m,n}$ that are the solutions of the following equations

$$\begin{aligned}
x_1 &= (g_{2,1} - g'_{1,1}) = g_{1,1}^2, \\
x_2 &= (g_{2,2} - g'_{1,2}) = g_{1,2}^2, \\
x_3 &= (g_{2,3} - g'_{1,3}) = g_{1,3}^2, \\
x_4 &= (g'_{2,1} + f_2 g_{2,1} + f_1 g_{1,1} + f_0) = -g_{1,1} g_{2,1}, \\
x_5 &= (g'_{2,2} + f_2 g_{2,2} + f_1 g_{1,2} + f_0) = -g_{1,2} g_{2,2}, \\
x_6 &= (g'_{2,3} + f_2 g_{2,3} + f_1 g_{1,3} + f_0) = -g_{1,3} g_{2,3},
\end{aligned} \qquad (48)$$

the system (42) takes the form

$$\begin{aligned}
y'_1 &= y_1 g_{1,1} - f \frac{g_{1,3} - g_{1,2}}{D} \\
y'_2 &= y_1 g_{1,2} - f \frac{g_{1,1} - g_{1,3}}{D} \\
y'_3 &= y_3 g_{1,3} - f \frac{g_{1,2} - g_{1,2}}{D}
\end{aligned} \qquad (49)$$

From (48) it follows that functions $g_{1,n}$ are the three different solutions of the system of the first-order nonlinear differential equations

$$\begin{aligned}
g'_{1,n} + g_{1,n}^2 - g_{2,n} &= 0 \\
g'_{2,n} + f_2 g_{2,n} + f_1 g_{1,n} + f_0 + g_{1,n} g_{2,n} &= 0
\end{aligned} \qquad (50)$$

This system can be written as the second-order nonlinear differential equation

$$g''_{1,n} + g'_{1,n}(3 g_{1,n} + f_2) + f_1 g_{1,n} + f_0 + f_2 g_{1,n}^2 + g_{1,n}^3 = 0. \qquad (51)$$

For the simplest case $f_2 \equiv 0, f_1 \equiv 0$

$$g''_{1,n} + g'_{1,n} 3 g_{1,n} + g_{1,n}^3 + f_0 = 0. \qquad (52)$$

In contrast to the Duffing equation (see, for example, [32]), the equation (52) has no linear term proportional to $g_{1,n}$.

For the third-order linear equation the equation (51) (or system (50)) plays the same role as the Riccati equation for second-order linear equation.

The functions $y_n = \exp\left(\int_{t_1}^{t} g_{1,n} dt'\right)$ are the linear independent ones, and the general solution of the homogeneous equation (36) ($f \equiv 0$) is

$$y(t) = \sum_{n=1}^{3} y_n(t_1) \exp\left(\int_{t_1}^{t} g_{1,n} dt'\right). \qquad (53)$$

There are other forms of system of the first order equations that can be obtained from the system (42) by choosing functions $g_{m,n}$.

## Conclusions

We presented some generalization of the already studied transformations of the linear differential equation into a system of the first order equations. The proposed transformation gives possibility to get new forms of the $N$-dimensional system of first order equations that can be useful for analysis of the solutions of the $N$-th order differential equation. In particular, for the third-order linear equation the nonlinear second-order equation that plays the same role as the Riccati equation for second-order linear equation is obtained.

**Appendix**

We represent the solution of the equation
$$y^{(N)} + f_{N-1}(t)y^{(N-1)} + \ldots f_1(t)y' + f_0(t)y + f(t) = 0. \tag{54}$$
as the sum of new functions
$$y(t) = \sum_{n=1}^{N} y_n(t). \tag{55}$$

By introducing $N$ new unknowns $y_n(t)$ instead of the one $y(t)$, we can impose additional conditions. These conditions we write in the form
$$\begin{aligned} y' &= \sum_{n=1}^{N} g_{1,n}(t) y_n(t) \\ y'' &= \sum_{n=1}^{N} g_{2,n}(t) y_n(t) \\ &\ldots \\ y^{(N-1)} &= \sum_{n=1}^{N} g_{N-1,n}(t) y_n(t) \end{aligned} \tag{56}$$

where $g_{m,n}(t)$ are the arbitrary continuous functions having continuous derivatives.

If
$$D(t) = \begin{vmatrix} 1 & 1 & .. & 1 \\ g_{1,1}(t) & g_{1,2}(t) & \ldots & g_{1,N}(t) \\ \ldots & \ldots & \ldots & \ldots \\ g_{N-1,1}(t) & g_{N-1,2}(t) & \ldots & g_{N-1,N}(t) \end{vmatrix} \neq 0, \tag{57}$$

then the representation (55)-(56) is unique. Indeed, from (55) and (56) we can uniquely find $y_n(t)$ as a linear combination if $y(t)$ and its derivatives. Equations (56) can be rewritten as
$$\begin{aligned} y' &= \sum_{n=1}^{N} y'_n(t) = \sum_{n=1}^{N} g_{1,n}(t) y_n(t), \\ y'' &= (y'_n)' = \sum_{n=1}^{N} g'_{1,n}(t) y_n(t) + \sum_{n=1}^{N} g_{1,n}(t) y'_n(t) = \sum_{n=1}^{N} g_{2,n}(t) y_n(t), \\ y''' &= (y'')' = \sum_{n=1}^{N} g'_{2,n}(t) y_n(t) + \sum_{n=1}^{N} g_{2,n}(t) y'_n(t) = \sum_{n=1}^{N} g_{3,n}(t) y_n(t), \\ &\ldots \\ y^{(N-2)} &= \left(y^{(N-3)}\right)' = \sum_{n=1}^{N} g'_{N-3,n}(t) y_n(t) + \sum_{n=1}^{N} g_{N-3,n}(t) y'_n(t) = \sum_{n=1}^{N} g_{N-2,n}(t) y_n(t), \\ y^{(N-1)} &= \left(y^{(N-2)}\right)' = \sum_{n=1}^{N} g'_{N-2,n}(t) y_n(t) + \sum_{n=1}^{N} g_{N-2,n}(t) y'_n(t) = \sum_{n=1}^{N} g_{N-1,n}(t) y_n(t). \end{aligned} \tag{58}$$

The final system of the first order difference equations has the form
$$\begin{aligned} \sum_{n=1}^{N} y'_n &= \sum_{n=1}^{N} g_{1,n} y_n, \\ \sum_{n=1}^{N} g_{m,n} y'_n &= \sum_{n=1}^{N} \left(g_{m+1,n} - g'_{m,n}\right) y_n, \quad 1 \leq m \leq N-2, \\ \sum_{n=1}^{N} g_{N-1,n} y'_n &= \sum_{n=1}^{N} L_n y_n - f(t), \end{aligned} \tag{59}$$

where
$$L_n = -g'_{N-1,n} - f_{N-1}g_{N-1,n} - f_{N-2}g_{N-2,n} - \ldots - f_0(t).  \quad (60)$$

In matrix form
$$Y'(t) = M^{-1}FY + M^{-1}H.  \quad (61)$$

where $Y = (y_1, y_2, \ldots, y_N)^T$, $H = (0, 0, \ldots, -f(t))^T$,

$$M = \begin{pmatrix} 1 & 1 & \ldots & 1 \\ g_{1,1} & g_{1,2} & \cdots & g_{1,N} \\ \ldots & \ldots & \ldots & .. \\ g_{N-2,1} & g_{N-2,2} & \ldots & g_{N-2,N} \\ g_{N-1,1} & g_{N-1,2} & \cdots & g_{N-1,N} \end{pmatrix}, \quad (62)$$

$$F = \begin{pmatrix} g_{1,1} & g_{1,2} & \cdots & g_{1,N} \\ g_{2,1} - g'_{1,1} & g_{2,2} - g'_{1,2} & \cdots & g_{2,N} - g'_{1,N} \\ \ldots & \ldots & \ldots & .. \\ g_{N-1,1} - g'_{N-2,1} & g_{N-1,2} - g'_{N-2,2} & \ldots & g_{N-1,N} - g'_{N-2,N} \\ L_1 & L_2 & \ldots & L_N \end{pmatrix}. \quad (63)$$